\pdfoutput=1
\RequirePackage{ifpdf}
\ifpdf 
\documentclass[pdftex]{sigma}
\else
\documentclass{sigma}
\fi

\usepackage[matrix,arrow,ps,color,line,curve,frame]{xy}

\begin{document}

\newcommand{\arXivNumber}{1402.3024}


\renewcommand{\thefootnote}{$\star$}

\renewcommand{\PaperNumber}{060}

\FirstPageHeading

\ShortArticleName{Local Proof of Algebraic Characterization of Free Actions}

\ArticleName{Local Proof of Algebraic Characterization\\
of Free Actions\footnote{This paper is a~contribution to the Special Issue on Noncommutative Geometry and Quantum Groups
in honor of Marc A.~Rief\/fel.
The full collection is available at
\href{http://www.emis.de/journals/SIGMA/Rieffel.html}{http://www.emis.de/journals/SIGMA/Rieffel.html}}}

\Author{Paul F.~BAUM~$^{\dag\ddag}$ and Piotr M.~HAJAC~$^{\ddag\S}$}

\AuthorNameForHeading{P.F.~Baum and P.M.~Hajac}

\Address{$^\dag$~Mathematics Department, McAllister Building, The Pennsylvania State University,\\
\hphantom{$^\dag$}~University Park, PA 16802, USA} 
\EmailD{\href{mailto:baum@math.psu.edu}{baum@math.psu.edu}}

\Address{$^\ddag$~Instytut Matematyczny, Polska Akademia Nauk, ul.~\'Sniadeckich 8, Warszawa, 00-656 Poland}
\EmailD{\href{mailto:pmh@impan.pl}{pmh@impan.pl}}
\URLaddressD{\url{http://www.impan.pl/~pmh}}

\Address{$^\S$~Katedra Metod Matematycznych Fizyki, Uniwersytet Warszawski,\\
\hphantom{$^\S$}~ul.~Ho\.za~74, 00-682 Warszawa, Poland}

\ArticleDates{Received February 13, 2014, in f\/inal form May 21, 2014; Published online June 06, 2014}

\Abstract{Let~$G$ be a~compact Hausdorf\/f topological group acting on a~compact Hausdorf\/f topological space~$X$.
Within the $C^{*}$-algebra $C(X)$ of all continuous complex-valued functions on~$X$, there is the Peter--Weyl algebra
$\mathcal{P}_G(X)$ which is the (purely algebraic) direct sum of the isotypical components for the action of~$G$
on~$C(X)$.
We prove that the action of~$G$ on~$X$ is free if and only if the canonical map
$\mathcal{P}_G(X)\otimes_{C(X/G)}\mathcal{P}_G(X)\to \mathcal{P}_G(X)\otimes\mathcal{O}(G)$ is bijective.
Here both tensor products are purely algebraic, and $\mathcal{O}(G)$ denotes the Hopf algebra of ``polynomial''
functions on~$G$.}


\Classification{22C05; 55R10; 57S05; 57S10}

\rightline{\it With admiration and affection, to Marc A.~Rieffel on the occasion of his 75th birthday}

\renewcommand{\thefootnote}{\arabic{footnote}} \setcounter{footnote}{0}

\section{Theorem}

Given a~compact Hausdorf\/f topological group~$G$, we denote by ${\mathcal O} (G)$ the dense Hopf $*$-subalgebra of the
commutative $C^{*}$-algebra $C(G)$ spanned by the matrix coef\/f\/icients of irreducible representations of~$G$.
Let~$X$ be a~compact Hausdorf\/f topological space with a~given continuous (right) action of~$G$.
The action map
\begin{gather*}
X\times G\ni (x,g)\longmapsto xg\in X
\end{gather*}
determines a~map of $C^{*}$-algebras
\begin{gather*}
\delta\colon \  C(X)\longrightarrow C(X\times G).
\end{gather*}
Moreover, denoting by $\otimes$ the purely algebraic tensor product over the f\/ield $\mathbb{C}$ of complex numbers, we
def\/ine the \emph{Peter--Weyl subalgebra}~\cite[(3.1.4)]{bhms} of $C(X)$ as
\begin{gather*}
\mathcal{P}_G(X):=\{a\in C(X)| \delta(a)\in C(X)\otimes{\mathcal O}(G)\}.
\end{gather*}
Using the coassociativity of $\delta$, one can check that $\mathcal{P}_G(X)$ is a~right ${\mathcal O} (G)$-comodule algebra.
In particular, $\mathcal{P}_G(G)={\mathcal O}(G)$.
The assignment $X\mapsto\mathcal{P}_G(X)$ is functorial with respect to continuous $G$-equivariant maps and comodule algebra homomorphisms.
We call it the \emph{Peter--Weyl functor}.
Equivalently~\cite[Proposition 2.2]{s-pm11}, $\mathcal{P}_G(X)$ is the (purely algebraic) direct sum of the isotypical
components for the action of~$G$ on~$C(X)$ (see~\cite[p.~31]{m-gd61} and~\cite{l-lh53}, cf.~\cite[Theorem~1.5.1]{p-p95}).
In the special case that the action of~$G$ on~$X$ is free, $\mathcal{P}_{G}(X)$ is the algebra of all continuous
sections of the vector bundle~$E$ on $X/G$, where
\begin{gather}\label{vecbun}
E:= X\underset{G}{\times}{\mathcal O}(G).
\end{gather}
Note that in forming this vector bundle, ${\mathcal O}(G)$ is topologized as the direct limit of its f\/inite-dimensional
vector subspaces, not by the norm topology.

The theorem of this paper is:

\begin{theorem}\label{theoGal}
Let~$X$ be a~compact Hausdorff space equipped with a~continuous $($right$)$ action of a~compact Hausdorff group~$G$.
Then the action is free if and only if the canonical map
\begin{gather}
{can} \colon \ \mathcal{P}_G(X)\underset{C(X/G)}{\otimes} \mathcal{P}_G(X) \longrightarrow \mathcal{P}_G(X)\otimes
\mathcal{O}(G),
\nonumber
\\
{can}\colon \ x\otimes y \longmapsto (x\otimes 1)\delta(y),
\label{pwg}
\end{gather}
is bijective.
$($Here both tensor products are purely algebraic.$)$
\end{theorem}
\begin{definition}
The action of a~compact Hausdorf\/f group~$G$ on a~compact Hausdorf\/f space~$X$ satisf\/ies the \emph{Peter--Weyl--Galois
$($PWG$)$ condition} if\/f the canonical map~\eqref{pwg} is bijective.
\end{definition}
\noindent Our result states that the usual formulation of free action is equivalent to the algebraic PWG-condition.
In particular, our result provides a~framework for extending Chern--Weil theory beyond the setting of dif\/ferentiable
manifolds and into the context of cyclic homology~--
noncommutative geometry~\cite{bh04}.

\section{Proof}

The proof of the equivalence of freeness and the PWG-condition consists of six steps.
The f\/irst step takes care of the easy implication of the equivalence, and the remaining f\/ive steps prove the dif\/f\/icult
implication of the equivalence.

\subsection[PWG-condition $\Rightarrow$ freeness]{PWG-condition $\boldsymbol{\Rightarrow}$ freeness}

It is immediate that the action is free, i.e.\ $xg=x \Longrightarrow g=e$ (where~$e$ is the identity element of~$G$), if
and only if
\begin{gather*}
F\colon \  X
\times
G  \longrightarrow X\underset{X/G}{
\times
}X,
\\
F\colon \  (x,g)  \longmapsto (x,xg),
\end{gather*}
is a~homeomorphism.
Here $X
\times
_{X/G}X$ is the subset of $X
\times
X$ consisting of pairs $(x_1, x_2)$ such that $x_1$ and $x_2$ are in the same~$G$-orbit.

This is equivalent to the assertion that the $*$-homomorphism
\begin{gather*}
F^{*}\colon \ C\Big(X\underset{X/G}{
\times
}X\Big)\longrightarrow C(X
\times
G)
\end{gather*}
obtained from the above map~$F$ is an isomorphism.
Note that~$F$ is always surjective, so that the $*$-homomorphism $F^{*}$ is always injective.
Furthermore, there is the following commutative diagram in which the vertical arrows are the evident maps:
\begin{gather}
\label{diagram1}
\begin{split}
& \xymatrix{
\mathcal{P}_{G}(X)\underset{C(X/G)}{\otimes}\mathcal{P}_{G}(X)\ar[r]^{\ \ \ {can}}\ar[d]^{}&
\mathcal{P}_{G}(X)\otimes{\mathcal O}(G)\ar[d]^{}
\\
C\Big(X\underset{X/G}{\times}X\Big)\ar[r]^{\ \ \ F^{*}}& C(X\times G).
}
\end{split}
\end{gather}

Since the right-hand side of the canonical map~\eqref{pwg}, i.e.\ $\mathcal{P}_G(X)\otimes \mathcal{O}(G)$, is dense in the $C^{*}$-algebra $C(X
\times
G)$, validity of the PWG-condition combined with the commutativity of the diagram~\eqref{diagram1} implies that the
image of the $*$-homomorphism $F^{*}$ is dense in $C(X
\times
G)$.
Therefore, as the image of a~$*$-homomorphism of $C^{*}$-algebras is always closed, PWG implies surjectivity of~$F^{*}$,
which in turn implies that the action of~$G$ on~$X$ is free.

\subsection{Reduction to surjectivity and matrix coef\/f\/icients}

\begin{lemma}\label{lemma}
Let~$X$ be a~compact Hausdorff space equipped with a~continuous $($right$)$ action of a~compact Hausdorff group~$G$.
Then the canonical map is surjective if and only if for any matrix coefficient~$h$ of an irreducible representation
of~$G$, the element $1\otimes h$ is in the image of the canonical map.
Moreover, if the canonical map is surjective, then it is bijective.
\end{lemma}

\begin{proof}
First observe that the canonical map is a~homomorphism of left $\mathcal{P}_{G}(X)$-modules.
The f\/irst assertion of the lemma follows by combining this observation with the fact that matrix coef\/f\/icients of
irreducible representations span ${\mathcal O}(G)$ as a~vector space.

The Hopf algebra ${\mathcal O}(G)$ is cosemisimple.
Hence, by the result of H.-J.~Schneider~\cite[Theorem~I]{s-hj90}, if the canonical map is surjective, then it is
bijective.
\end{proof}

Alternately, assuming that the action of~$G$ on~$X$ is free, injectivity can be directly proved by using the
vector-bundle point of view indicated in~\eqref{vecbun} (see~\cite{bdh}).

\subsection{Reduction to free actions of compact Lie-groups}
\label{reduction}

Assume that Theorem~\ref{theoGal} holds for compact Lie groups.
In this section, we prove that this special case implies Theorem~\ref{theoGal} in general.

Let $\varphi\colon G\to U(n)$ be any f\/inite-dimensional representation of~$G$.
Set
\begin{gather*}
X_\varphi:=X\underset{G}{\times }U(n).
\end{gather*}
Thus $X_\varphi= (X
\times
U(n))/G$, where~$G$ acts on $X
\times
U(n)$ by $(x,u)g:= (xg,\varphi(g^{-1})u) $.
The group $U(n)$ acts on $X_\varphi$ by $[(x,u)]v:=[(x,uv)]$, and this action is free.
The map
\begin{gather*}
\Phi\colon \ X\longrightarrow X_\varphi,
\qquad
x\longmapsto [(x,I_n)],
\end{gather*}
where $I_n\in U(n)$ is the identity matrix, has the equivariance property
\begin{gather}
\label{equivariance}
\Phi(xg)=\Phi(x)\varphi(g).
\end{gather}
Hence $\Phi$ and $\varphi$ induce maps $\Phi^{*}\colon \mathcal{P}_{U(n)}(X_\varphi)\to \mathcal{P}_{G}(X)$ and
$\varphi^{*}\colon{\mathcal O}(U(n))\to{\mathcal O}(G)$.
The equivariance property~\eqref{equivariance} implies commutativity of the diagram
\begin{gather*}
\xymatrix{
\mathcal{P}_{U(n)}(X_\varphi)\underset{C(X/G)}{\otimes}\mathcal{P}_{U(n)}(X_\varphi)
\ar[r]^{\ \ \ {can}}\ar[d]^{\Phi^{*}\otimes\Phi^{*}}&
\mathcal{P}_{U(n)}(X_\varphi)\otimes{\mathcal O}(U(n))\ar[d]^{\Phi^{*}\otimes\varphi^{*}}\\
\mathcal{P}_{G}(X)\underset{C(X/G)}{\otimes}\mathcal{P}_{G}(X)\ar[r]^{\ \ \ {can}}& \mathcal{P}_{G}(X)\otimes{\mathcal O}(G).
}
\end{gather*}
Therefore surjectivity of the upper canonical map implies that $1\otimes h$ is in the image of the lower canonical map,
where~$h$ is any matrix coef\/f\/icient of~$\varphi$.
By Lemma~\ref{lemma}, this implies the PWG-condition.

\subsection{Local triviality for free actions of compact Lie groups}

We recall the theorem of A.M.~Gleason:
\begin{theorem}[\cite{g-am50}] Let~$G$ be a~compact Lie group acting freely and continuously on a~completely regular space~$X$.
Then~$X$ is a~locally trivial~$G$-bundle over $X/G$.
\end{theorem}

Combining the Gleason theorem with Section~\ref{reduction}, we infer that the PWG-condition is valid for free actions if it
is valid for locally trivial free actions.

\subsection{Reduction to the trivial-bundle case}\label{tbundle}

Assume that the action of~$G$ on~$X$ is free and locally trivial.
Since the quotient space $X/G$ is compact Hausdorf\/f, we can choose a~f\/inite open cover $U_1,\ldots,U_r$ of $X/G$ such
that each $\pi^{-1}(U_j)$ is a~trivializable principal~$G$-bundle over~$U_j$.
Here $\pi\colon X\to X/G$ is the quotient map.
On $X/G$, let $\psi_1,\ldots,\psi_r$ be a~partition of unity subordinate to the cover $U_1,\ldots,U_r$.
Then, for each~$j$ there is the canonical map
\begin{gather*}
{can}_j\colon \mathcal{P}_G\big(\pi^{-1}(\mathrm{supp}(\psi_j)\big)\underset{C(\mathrm{supp}(\psi_j))} {\otimes}
\mathcal{P}_G\big(\pi^{-1}(\mathrm{supp}(\psi_j)\big) \longrightarrow
\mathcal{P}_G\big(\pi^{-1}(\mathrm{supp}(\psi_j)\big)\otimes{\mathcal O}(G).
\end{gather*}
Assume that for each $j\in\{1,\ldots,r\}$ there exist elements
\begin{gather*}
p_{j1},q_{j1},\ldots,p_{jn},q_{jn}\in \mathcal{P}_G\big(\pi^{-1}(\mathrm{supp}(\psi_j)\big)
\qquad
\text{such that}
\qquad
{can}_j\left(\sum\limits_{i=1}^n p_{ji}\otimes q_{ji}\right)=1\otimes h.
\end{gather*}
Let $\widetilde{p_{ji}}$'s and $\widetilde{q_{ji}}$'s be functions on~$X$ obtained respectively from functions
$p_{ji}$'s and $q_{ji}$'s by zero-value extension.
Then for each~$i$ and~$j$ we take
\begin{gather*}
\widetilde{p_{ji}}\sqrt{\psi_j\circ\pi},\; \widetilde{q_{ji}}\sqrt{\psi_j\circ\pi}
\in\mathcal{P}_G(X),
\end{gather*}
and for any $x\in X$ and $g\in G$, using the commutativity of the diagram~\eqref{diagram1}, we obtain
\begin{gather*}
{can}  \left(\sum\limits_{j=1}^r\sum\limits_{i=1}^n \widetilde{p_{ji}}\sqrt{\psi_j\circ\pi}\otimes
\widetilde{q_{ji}}\sqrt{\psi_j\circ\pi} \right) (x,g)\\
\qquad{}
 = \left(\sum\limits_{j=1}^r\sum\limits_{i=1}^n \widetilde{p_{ji}}\sqrt{\psi_j\circ\pi}\otimes
\widetilde{q_{ji}}\sqrt{\psi_j\circ\pi} \right) (x,xg)
\\
\qquad{}
 =  \sum\limits_{\text{all}\; j\;\text{s.t.}\; \pi(x)\in U_j} (\psi_j\circ\pi)(x)
\sum\limits_{i=1}^n p_{ji}(x) q_{ji}(xg)
\\
\qquad{}
 =  \sum\limits_{\text{all}\; j\;\text{s.t.}\; \pi(x)\in U_j} (\psi_j\circ\pi)(x) \, {can}_j
\left(\sum\limits_{i=1}^n p_{ji}\otimes q_{ji}\right)(x,g)
\\
\qquad{}
 = \sum\limits_{j=1}^r (\psi_j\circ\pi)(x) h(g)
\\
\qquad{}
 =(1\otimes h)(x,g).
\end{gather*}
Hence validity of the PWG-condition in the trivial-bundle case implies that the PWG-condition holds for actions that are
free and locally trivial.

\subsection{The trivial-bundle case}

First note that
\begin{gather}\label{note}
\mathcal{P}_{G}(Y\times G)=C(Y)\otimes\mathcal{O}(G).
\end{gather}
This is implied by two facts: (1) quite generally $\mathcal{P}_G(X)\subseteq C(X)$ is the purely algebraic direct sum of
the isotypical components for the action of G on $C(X)$~\cite[Proposition 2.2]{s-pm11}; and (2)~each isotypical
component for the action of G on $C(Y
\times
G)$ is of the form $C(Y)\otimes V$, where~$V$ is an isotypical component for the action of~$G$ on~$C(G)$.

As ${\mathcal O}(G)$ is a~Hopf algebra, the dual of the homeomorphism
\begin{gather*}
F_G\colon \  G\times G\ni (g_1,g_2)\longmapsto (g_1,g_1g_2)\in G\times G
\end{gather*}
and the dual of its inverse $F_G^{-1}$ restrict and corestrict respectively to
\begin{gather*}
{can}_G\colon \ {\mathcal O}(G)\otimes{\mathcal O}(G)\ni T \longmapsto T\circ F_G\in {\mathcal O}(G)\otimes{\mathcal O}(G),
\\
{can}_G^{-1}\colon \  {\mathcal O}(G)\otimes{\mathcal O}(G)\ni T \longmapsto T\circ F_G^{-1}\in {\mathcal
O}(G)\otimes{\mathcal O}(G).
\end{gather*}
Granted the identif\/ication~\eqref{note}, we now obtain the following commutative diagram:
\begin{gather*}
\xymatrix{
\mathcal{P}_{G}(Y\times G)\underset{C(Y)}{\otimes}\mathcal{P}_{G}(Y\times G) \ar[r]^{\qquad {can}}\ar[d]^{=}&
\mathcal{P}_{G}(Y\times G)\otimes{\mathcal O}(G)\ar[d]^{=}\\
C(Y)\otimes{\mathcal O}(G)\otimes{\mathcal O}(G)\ar[r]^{\mathrm{id} \otimes {can}_G}& C(Y)\otimes{\mathcal O}(G)\otimes{\mathcal O}(G).
}
\end{gather*}
Hence the bijectivity of ${can}_G$ implies the bijectivity of~${can}$.

\section{Appendix}

In this appendix, we observe that Step~\ref{tbundle}, i.e.\ reduction to the trivial-bundle case, is implied by the
following general results:
\begin{lemma}
Let $(H,\Delta)$ be a~compact quantum group acting on a~unital $C^{*}$-algebra~$A$.
The assignment $A\mapsto \mathcal{P}_H(A)$ of the Peter--Weyl algebra to a~$C^{*}$-algebra yields a~functor from the
category with objects being unital $C^{*}$-algebras with an~$H$-coaction $($and morphisms being equivariant unital
$^{*}$-homomorphisms$)$ to the category whose objects are ${\mathcal O}(H)$-comodule algebras $($and whose morphisms are
colinear algebra homomorphisms$)$.
Furthermore, this functor commutes with all $($equivariant$)$ pullbacks, that is,
\begin{gather*}
\mathcal{P}_H\Big(A\underset{B}{\times} C\Big)=\mathcal{P}_H(A)\underset{\mathcal{P}_H(B)}{\times}\mathcal{P}_H(C).
\end{gather*}
\end{lemma}

\begin{lemma}[Lemma~3.2 in~\cite{hkmz11}]
Let $\mathcal{H}$ be a~Hopf algebra with bijective antipode, and let
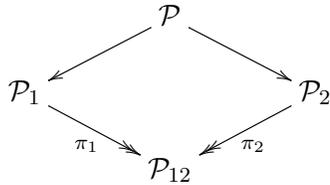
\begin{gather*}
\mbox{
$\xymatrix@=5mm{& & \mathcal{P} \ar[lld]
\ar[rrd] & &\\
\mathcal{P}_1 \ar@{>>}[drr]_{\pi_1}& & & &\mathcal{P}_2 \ar@{>>}[dll]^{\pi_2}\\
&& \mathcal{P}_{12} &&}$
}
\end{gather*}
be the pullback diagram of surjective right $\mathcal{H}$-comodule algebra homomorphisms.
Then $\mathcal{P}$ is principal if and only if $\mathcal{P}_1$ and $\mathcal{P}_2$ are principal\footnote{For the def\/inition of ``principal'' see~\cite{bh04} and~\cite[Def\/inition~2.3]{hkmz11}.}.
\end{lemma}

Proving the f\/irst lemma is straightforward, and the second lemma is the highlight of~\cite{hkmz11}.
For the case considered in this paper (a compact Hausdorf\/f group~$G$ acting continuously on a~compact Hausdorf\/f
space~$X$), we have $A=C(X)$, $H=C(G)$, ${\mathcal O}(H)={\mathcal O}(G)$, $\mathcal{P}_H(A)=\mathcal{P}_G(X)$, and the
condition of being ``principal'' is equivalent to the PWG-condition.
Thus these two lemmas combined with standard induction yield an alternative proof of Step~\ref{tbundle}.

Moreover, the theorem of this paper is a~special case of a~much more general theorem about compact quantum groups acting
on unital $C^{*}$-algebras~\cite{bdh}.
However, the proof of the general theorem is nonlocal.

\subsection*{Acknowledgments}

We thank the referees for the careful attention they have given to this paper.
This work was partially supported by NCN grant 2011/01/B/ST1/06474.
P.F.~Baum was partially supported by NSF grant DMS 0701184.

\pdfbookmark[1]{References}{ref}
\LastPageEnding


\begin{thebibliography}{99}
\footnotesize\itemsep=0pt

\bibitem{bdh}
Baum P.F., De~Commer K., Hajac P.M., Free actions of compact quantum groups on
  unital $C^{*}$-algebras, \href{http://arxiv.org/abs/1304.2812}{arXiv:1304.2812}.

\bibitem{bhms}
Baum P.F., Hajac P.M., Matthes R., Szyma\'nski W., Noncommutative geometry
  approach to principal and associated bundles, in Quantum Symmetry in
  Noncommutative Geometry, {t}o appear, \href{http://arxiv.org/abs/math.DG/0701033}{math.DG/0701033}.

\bibitem{bh04}
Brzezi{\'n}ski T., Hajac P.M., The {C}hern--{G}alois character,
  \href{http://dx.doi.org/10.1016/j.crma.2003.11.009}{\textit{C.~R.~Math. Acad. Sci. Paris}} \textbf{338} (2004), 113--116,
  \href{http://arxiv.org/abs/math.KT/0306436}{math.KT/0306436}.

\bibitem{g-am50}
Gleason A.M., Spaces with a compact {L}ie group of transformations,
  \href{http://dx.doi.org/10.1090/S0002-9939-1950-0033830-7}{\textit{Proc. Amer. Math. Soc.}} \textbf{1} (1950), 35--43.

\bibitem{hkmz11}
Hajac P.M., Kr{\"a}hmer U., Matthes R., Zieli{\'n}ski B., Piecewise principal
  comodule algebras, \href{http://dx.doi.org/10.4171/JNCG/88}{\textit{J.~Noncommut. Geom.}} \textbf{5} (2011), 591--614,
  \href{http://arxiv.org/abs/0707.1344}{arXiv:0707.1344}.

\bibitem{l-lh53}
Loomis L.H., An introduction to abstract harmonic analysis, D.~Van Nostrand
  Company, Inc., Toronto~-- New York~-- London, 1953.

\bibitem{m-gd61}
Mostow G.D., Cohomology of topological groups and solvmanifolds, \href{http://dx.doi.org/10.2307/1970281}{\textit{Ann.
  of Math.}} \textbf{73} (1961), 20--48.

\bibitem{p-p95}
Podle{\'s} P., Symmetries of quantum spaces. {S}ubgroups and quotient spaces of
  quantum {${\rm SU}(2)$} and {${\rm SO}(3)$} groups, \href{http://dx.doi.org/10.1007/BF02099436}{\textit{Comm. Math.
  Phys.}} \textbf{170} (1995), 1--20, \href{http://arxiv.org/abs/hep-th/9402069}{hep-th/9402069}.

\bibitem{s-hj90}
Schneider H.-J., Principal homogeneous spaces for arbitrary {H}opf algebras,
  \href{http://dx.doi.org/10.1007/BF02764619}{\textit{Israel~J. Math.}} \textbf{72} (1990), 167--195.

\bibitem{s-pm11}
So{\l}tan P.M., On actions of compact quantum groups, \textit{Illinois~J.
  Math.} \textbf{55} (2011), 953--962, \href{http://arxiv.org/abs/1003.5526}{arXiv:1003.5526}.

\end{thebibliography}
\end{document}